\documentclass{amsart}

\usepackage{amssymb}
\usepackage[all]{xy}

\newtheorem{thm}{Theorem}%[section]

\newtheorem{lem}[thm]{Lemma}
\newtheorem{cor}[thm]{Corollary}

\newtheorem{prop}[thm]{Proposition}

\theoremstyle{definition}
\newtheorem{defn}[thm]{Definition}

\newtheorem{say}[thm]{}

\newtheorem{ques}[thm]{Question}    %!!!!!!!!!!!!!!!!!!!!

\newtheorem{rem}[thm]{Remark}          

\newtheorem*{ack}{Acknowledgments}      % \renewcommand{\theack}{} 

\newtheorem{defn-thm}[thm]{Definition--Theorem}  %!!!!!!!!!!!!!!!!!!!!!!!!
\newtheorem{defn-lem}[thm]{Definition--Lemma}  %!!!!!!!!!!!!!!!!!!!!!!!!
\theoremstyle{remark}

\renewcommand{\c}[0]{{\mathbb C}}  

\renewcommand{\o}[0]{{\mathcal O}} 

\newcommand{\n}[0]{{\mathbb N}}
\newcommand{\dd}[0]{{\mathbb D}}

\newcommand{\qtq}[1]{\quad\mbox{#1}\quad}

\newcommand{\mult}[0]{\operatorname{mult}}

\newcommand{\ex}[0]{\operatorname{Ex}}

\newcommand{\der}[0]{\operatorname{Der}}

\newcommand{\tsum}[0]{\textstyle{\sum}}

\def\loccoh#1.#2.#3.#4.{H^{#1}_{#2}(#3,#4)}

\DeclareMathAlphabet{\mathchanc}{OT1}{pzc}%
                                {m}{it}

\usepackage[all]{xy}\xyoption{dvips}

\begin{document}
\bibliographystyle{amsalpha}

%\today

 \title{Partial resolution by toroidal blow-ups}
 \author{J\'anos Koll\'ar}
\begin{abstract} We  give an alternate proof of a theorem of Tevelev about improving a non-toroidal ideal sheaf by a sequence of toroidal blow-ups.
\end{abstract}

 \maketitle

\begin{say}[Toroidal blow-up] \label{torbu.strat.defn}
Let $X$ be a smooth variety over a field  and $\sum D_i$ a
simple normal crossing  (abbreviated as {\it snc}) divisor on $X$.
A (closed) {\it stratum} of $(X, \sum D_i)$ is an 
irreducible component of an intersection $D_{i_1}\cap\cdots\cap D_{i_r}$. 
If $Z\subset X$ is a  stratum (or a disjoint union of strata)
and $\pi:B_ZX\to X$ the blow-up
then $\bigl(B_ZX, \sum_i \pi^{-1}_*D_i+\sum_j E_j\bigr)$
is also an snc pair where the $E_j$ are the  exceptional divisors of $\pi$. 
We call such blow-ups {\it toroidal.}
\end{say}

The following question was suggested by Keel.

\begin{ques}  \label{ques.1}
Let $(X, \sum D_i)$ be an snc pair over a field and $J\subset \o_X$ an ideal sheaf. How much can one improve $J$ by a sequence of toroidal blow-ups?
\end{ques}

As a simple example, assume that $X$ is a surface.  Then there are very few toroidal blow-ups:
we can  blow up either the curves $D_i\subset X$
(giving the identity map) or any of their intersection points. 
Thus if the cosupport of $J$ (that is, the support of $\o_X/J$)  does not contain any  strata then
toroidal blow-ups have no effect on $J$. 
Similarly, one expects to be able to improve the singularities of $J$ along strata but not necessarily along other subvarieties.  This leads to the following.

\begin{defn}  Let $(X, \Delta:=\sum D_i)$ be an snc pair over a field  and $J\subset \o_X$ an ideal sheaf. We say that $J$  is {\it toroidally resolved} if its cosupport does not contain any strata.
\end{defn}

The key step of the proof is to show that each ideal sheaf $J\subset \o_X$ has a
unique {\it toroidal hull} $J\subset J^t\subset \o_X$ such that the 
toroidal resolution problem for $J$ is equivalent to the
ordinary resolution problem for $J^t$; see Definition \ref{TOR.HULL.EXISTS}
  and 
 Proposition \ref{main.equiv.prop}.
The resolution of  toroidal ideals is known over arbitrary fields by
\cite{MR2219845}, thus we get the 
 following  answer to Question \ref{ques.1}.

\begin{thm} \label{main.thm.1}
Let $(X, \Delta)$ be an snc pair over a field (of arbitrary characteristic)  and $J\subset \o_X$ an ideal sheaf.
Then there is a toroidal blow-up sequence
$$
\bigl(X_n, \Delta_n, J_n\bigr)\to \cdots  \to \bigl(X_0,  \Delta_0, J_0\bigr):=
\bigl(X, \Delta, J\bigr)
$$
such that   $J_n\subset \o_{X_n}$ is  toroidally   resolved.
\end{thm}

We  state a more precise version in  Theorem \ref{main.thm.2}
and also explain how the ideals $J_i$ transform into each other, but 
first we
apply  Theorem \ref{main.thm.1} to the ideal sheaf of a divisor to get the following answer to the
original  question of Keel.

Tevelev pointed out that, using \cite{MR803344},  the methods of 
\cite{MR2343384} can easily be modified to  obtain
Corollary \ref{main.thm.1.cor};
 see also \cite{MR2452307, MR3343903, 2015arXiv151102203V} for
closely related variants.
In fact, \cite{MR2343384} gives the stronger result that 
 $\Pi^{-1}_*Y$ intersects each stratum in the expected codimension. 

\begin{cor}[Tevelev] \label{main.thm.1.cor}
Let $(X, \Delta)$ be an snc pair over a field  and $Y\subset X$ a closed subscheme that does not contain any of the
irreducible components of $\Delta$.
Then there is a sequence of toroidal blow-ups 
$\Pi: X_n\to \cdots  \to X_0:= X$
such that  the birational transform  $\Pi^{-1}_*Y$ does not contain any strata of  the pair $\bigl(X_n, \Pi^{-1}_*\Delta+\ex(\Pi)\bigr)$. \qed
\end{cor}

For another application, note that if a divisor  $B$  does not contain any strata of $(X, \Delta)$ 
iff $(X, \Delta+\epsilon B)$ is divisorial log terminal  (abbreviated as {\it dlt}) for $0<\epsilon\ll 1$, cf.\ \cite[2.8]{kk-singbook}. 
We can thus restate the divisorial case of Corollary \ref{main.thm.1.cor} as follows.

\begin{cor} \label{main.thm.1.cor.2}
Let $(X, \Delta)$ be an snc pair over a field  and $B\subset X$ an effective  divisor that does not contain any of the
irreducible components of $\Delta$.
Then there is a sequence of toroidal blow-ups 
$\Pi: X_n\to \cdots  \to X_0:= X$
such that  
$$
\bigl(X_n, \Pi^{-1}_*(\Delta+\epsilon B)+\ex(\Pi)\bigr)
\qtq{is  dlt   for} 0<\epsilon\ll 1. \qed
$$
\end{cor} 

The model obtained in  Corollary \ref{main.thm.1.cor.2} is related to the 
dlt modifications of $(X, \Delta+\epsilon B)$ constructed in \cite{MR2955764}
(in characteristic 0). Our models are smooth but the log canonical class need not be relatively nef. Nonetheless, 
this suggests that Corollary \ref{main.thm.1.cor.2} might be approached using the minimal model program. A problem is that there are many different 
 dlt modifications and most of them are singular. It is not clear to me how to guarantee smoothness using MMP.

\begin{say}[Plan of the proof of Theorem \ref{main.thm.1}] Assume for simplicity that
$(X, \Delta)$ is toric with torus $T$. We assume that $\Delta$ consists of all $T$-invariant divisors. We show that Theorem \ref{main.thm.1} for $J$ is essentially equivalent to a special case of  resolution, usually called monomialization, of the toric ideal
$J^t:=\sum_{\tau}\tau^*J$ where we sum over all $\tau\in T$. The latter is a combinatorial problem that is independent of the characteristic.

In general, $(X, \Delta)$ is locally toric  in the analytic or \'etale topology so we need to check that the  local construction of $J^t$
gives a global ideal sheaf $J^t$. This is probably well  known to experts.
I do not know a reference that covers everything that we need, so we go through the details.
\end{say}

In the  precise version of Theorem \ref{main.thm.1}  we further restrict the 
blow-ups allowed in the sequence. For this we need some definitions first.

\begin{say}[Toroidally equimultiple blow-ups] \label{equim.bu.say}
Let $X$ be a smooth variety and $J\subset \o_X$ an ideal sheaf.
Let $Z\subset X$ be a smooth subvariety and
$\pi:B_ZX\to X$ the blow-up of $Z$. Let $E\subset B_ZX$ denote the exceptional divisor.

Most resolution methods work with blow-up 
 centers $Z\subset X$ such that  $J$ is equimultiple along $Z$; that is,
$\mult_zJ=m$ for every $z\in Z$ for some fixed $m$. We then define the
{\it birational transform} of $J$ by
$$
\pi^{-1}_*J:=\o_{B_ZX}(mE)\cdot \pi^*J.
\eqno{(\ref{equim.bu.say}.1)}
$$
(This is frequently called the `controlled'  or `weak' transform.)
This is an ideal sheaf on $B_ZX$. It has the pleasant property
that $\mult_y\pi^{-1}_*J\leq m$ for every $y\in E$. 

Working toroidally,  we would like $Z$ to be a stratum (or a disjoint union of strata). However, if the multiplicity of $J$ jumps at a single point that is not a stratum, then  toroidal blow-ups are unlikely to  change this. Thus, in a resolution procedure,  the best one can hope for is that
 $J$ is   {\it toroidally equimultiple} along $Z$, that is,
$\mult_WJ=\mult_ZJ$ for every stratum $W\subset Z$.  

If this holds then we define the
{\it birational transform} of $J$ by
$$
\pi^{-1}_*J:=\o_{B_ZX}(mE)\cdot \pi^*J.
\eqno{(\ref{equim.bu.say}.2)}
$$
As before, this is an ideal sheaf on $B_ZX$ and
$\mult_V\pi^{-1}_*J\leq m$ for every stratum $V\subset E$. 

The resulting birational transform of $J$ 
then behaves as expected over generic points of strata $W\subset Z$ but 
can be rather badly behaved elsewhere. This is not a problem if we care only about generic points of strata.
\end{say}

Let us recall a somewhat detailed form of
resolution (usually called monomialization)  of ideal sheaves.

\begin{thm} \cite[3.68]{k-res}\label{main.thm.class}
Let $(X, E)$ be an snc pair over a field of characteristic 0 and $J\subset \o_X$ an ideal sheaf.
Then there is a  blow-up sequence
$$
\bigl(X_n, J_n, E_n\bigr)\to \cdots  \to \bigl(X_0, J_0, E_0\bigr):=
\bigl(X, J, E\bigr)
$$
with the following properties.  
\begin{enumerate}
\item Each $\pi_i:X_{i+1}\to X_i$ is   a blow-up with smooth center
$Z_i\subset X_i$ and exceptional divisor $E^{i+1}$.
\item $J_i$ is  equimultiple along $Z_i$.
\item $J_{i+1}=(\pi_i)^{-1}_*J_i$ as in (\ref{equim.bu.say}.1).
\item $Z_i$ has normal crossings with $E_i$ and $E_{i+1}=(\pi_i)^{-1}_*E_i+E^{i+1}$.
\item $\bigl(X_n, J_n, E_n\bigr) $ is    resolved; that is,
$J_n=\o_{X_n}$.  
\end{enumerate}
\end{thm}

Now we can state the more precise form of Theorem \ref{main.thm.1}
where we just add `toroidal' to the formulation of 
Theorem \ref{main.thm.class} in a  few places.

\begin{thm} \label{main.thm.2}
Let $(X, \Delta)$ be an snc pair over a field of any characteristic  and $J\subset \o_X$ an ideal sheaf.
Then there is a toroidal blow-up sequence
$$
\bigl(X_n,  \Delta_n, J_n\bigr)\to \cdots  \to \bigl(X_0, \Delta_0, J_0\bigr):=
\bigl(X, \Delta, J\bigr)
$$
with the following properties.  
\begin{enumerate}
\item Each $\pi_i:X_{i+1}\to X_i$ is   a blow-up with smooth, toroidal center
$Z_i\subset X_i$ and exceptional divisor $E_{i+1}$.
\item $J_i$ is toroidally equimultiple along $Z_i$.
\item $J_{i+1}=(\pi_i)^{-1}_*J_i$ as in (\ref{equim.bu.say}.2).
\item $\Delta_{i+1}=(\pi_i)^{-1}_*\Delta_i+E_{i+1}$.
\item $\bigl(X_n, \Delta_n, J_n\bigr) $ is toroidally   resolved.
\end{enumerate}
\end{thm} 

\begin{rem} The role of the divisors $E$ and $\Delta$ is quite different in the two Theorems; the notation is changed to emphasize this. 
In Theorem \ref{main.thm.class} $E$ is but an auxiliary datum which gives very mild restrictions on the blow-up centers, whereas 
in Theorem \ref{main.thm.1}  $\Delta$ gives extremely strong 
restrictions on the blow-up centers. 
\end{rem}

\begin{defn} Let us call a blow-up sequence 
satisfying (\ref{main.thm.class}.1--4) 
{\it equimultiple} and
 a blow-up sequence 
satisfying (\ref{main.thm.2}.1--4) 
{\it toroidally equimultiple.}

Thus Theorem \ref{main.thm.class} says that, in characteristic 0,  every ideal sheaf can be resolved by an equimultiple
  blow-up sequence.
\end{defn}

\begin{say}[Toroidal ideals]\label{tor.ideal.say}
Let $X$ be a smooth variety and $\sum D_i$ an snc divisor. An ideal sheaf
$I\subset \o_X$ is {\it toroidal} if $X$ is covered by open sets $U_j$
such that 
$$
I|_{U_j}=\sum_s \o_{U_j}\bigl(-\tsum_i m_{ijs}D_i|_{U_j}\bigr)
\eqno{(\ref{tor.ideal.say}.1)}
$$
for every $j$ and for suitable  $m_{ijs}\in \n$.

Let $Z\subset X$ be a closed stratum and 
$Z^0:=Z\setminus\cup\{W: W\subsetneq Z\mbox{ is a stratum}\}$ 
the corresponding {\it open stratum.}
For every $z\in Z^0\cap U_j$ the $m_{ijs} $ give vectors
$$
v_{js}:=\bigl(m_{ijs}: D_i\supset Z\bigr)\in \tsum_{i: D_i\supset Z}\n[D_i]
\eqno{(\ref{tor.ideal.say}.2)}
$$
and these generate a subsemigroup 
$$
M_Z\subset  \tsum_{i: D_i\supset Z}\n[D_i]
\eqno{(\ref{tor.ideal.say}.3)}
$$
which depends only on $Z$. For any inclusion of strata $W\subset Z$
we have the  coordinate  projection
$$
p_{Z,W}: \tsum_{i: D_i\supset W}\n[D_i]\to \tsum_{i: D_i\supset Z}\n[D_i]
\eqno{(\ref{tor.ideal.say}.4)}
$$
and the subsemigoups $M_Z$ satisfy the compatibility relation
$$
p_{Z,W}\bigl(M_W\bigr)=M_Z.
\eqno{(\ref{tor.ideal.say}.5)}
$$
This gives a one-to-one correspondence between toroidal ideals and
collections of subsemigroups  $\{M_Z\}$ satisfying the
compatibility relations
(\ref{tor.ideal.say}.5).
In particular, we see that $I\mapsto I^{\rm an}$ gives a one-to-one correspondence
$$
\{\mbox{toroidal ideals $I\subset \o_X$}\}
\ \leftrightarrow\
\{\mbox{toroidal ideals $I^{\rm an}\subset \o^{\rm an}_X$}\}.
\eqno{(\ref{tor.ideal.say}.6)}
$$
\end{say}

We claim that toroidal ideals are the only ones that can be
`canonically' associated to the stratification of an snc pair.

\begin{say}[Local stratified isomorphisms] \label{lsi.say}
Let $(X, \Delta)$ be  an snc pair and $U_1, U_2\subset X$   open sets. An isomorphism
$\phi:U_1\to U_2$ is called {\it stratification preserving} if
$Z\cap U_1=\phi^{-1}(Z\cap U_2)$ for every stratum $Z\subset X$. 
Note that our strata are the irreducible components of the intersections of the $D_i$, thus this is stronger than just assuming 
$D_i\cap U_1=\phi^{-1}(D_i\cap U_2)$ for every  $D_i$. 

We say that an ideal sheaf $I\subset \o_X$ is invariant under
stratification preserving local isomorphisms if
$\phi^*\bigl(I|_{U_2}\bigr)=I|_{U_1}$ holds for every such $\phi:U_1\to U_2$.
 
It is clear that a toroidal ideal is invariant under
stratification preserving local isomorphisms and
we would like to claim the converse.
Unfortunately, if $X$ has no birational automorphisms then the identity map is the only 
stratification preserving local isomorphism. As usual, there are 3 ways to get more $U_i$.
\medskip

{\it Complex analytic \ref{lsi.say}.1.} If $X$ is over $\c$, we  use
analytic open sets $U_1, U_2\subset X^{\rm an}$.

\medskip

{\it Etale local \ref{lsi.say}.2.}  We use \'etale morphisms
$\tau_i:U\to X$ and require that $\tau_1^{-1}(Z)=\tau_2^{-1}(Z)$
for every stratum $Z\subset X$. 
\medskip

{\it Formal local \ref{lsi.say}.3.} We use isomorphisms of complete local rings
$\phi^*:\hat{\o}_{x_2, X}\to \hat{\o}_{x_1, X}$. (If the base field is not algebraically closed we also allow residue field extensions.)

\medskip

{\it Micro local \ref{lsi.say}.4.} We assume the condition on the tangent space level. That is
$$
\der_X\bigl(-\log \Delta\bigr)\cdot I\subset I
$$
where $\der_X\bigl(-\log \Delta\bigr) $ is the sheaf of 
logarithmic derivatives 
 along $\Delta$;  cf.\ \cite[3.87]{kk-singbook}. This works in characteristic 0 but not in positive characteristic.
This shows that the concepts of toroidal ideal and
toroidal hull  (\ref{TOR.HULL.EXISTS}) are related to
D-balanced ideals and  well-tuned  ideals used in resolution.
See \cite[Sec.3.4]{k-res} for the latter notions.
\end{say}

\begin{prop} \label{cam.=.tor.prop}
Let $(X, \Delta)$ be  an snc pair and 
$I\subset \o_X$ an ideal sheaf that is invariant under all
stratification preserving local isomorphisms in any of the settings
(\ref{lsi.say}.1--3). Then $I$ is a toroidal ideal sheaf.
\end{prop}

Proof. We explain the complex analytic case and leave the details of the other settings to the reader. 
By (\ref{tor.ideal.say}.6) it is enough to show that $I^{\rm an}$ is toroidal.

Let $\dd\subset \c$ denote the unit disc and $\dd^*$ the punctured unit disc.
We will view $\dd^*\subset \c^*$ as a semigroup.

Let $Z^0\subset X$ be an open stratum.  
After reindexing the $D_i$, for every $z\in Z^0$ we can choose a
neighborhood of the form $(0\in \dd^n)$ where $D_i=(x_i=0)$ for $i=1,\dots, m$. 
We start with the natural   $(\dd^*)^m$ action on the first $m$ coordinates.
This is a  stratification preserving action.

Pick  any $f=\sum_{i_1,\dots, i_m}f_{i_1,\dots, i_m}(x_{m+1},\dots, x_n)\cdot x_1^{i_1}\cdots x_m^{i_m}\in I^{\rm an}$.
Then 
$$
\tau^*f=\sum_{i_1,\dots, i_m}\chi_{i_1,\dots, i_m}\cdot f_{i_1,\dots, i_m}(x_{m+1},\dots, x_n)\cdot x_1^{i_1}\cdots x_m^{i_m}
$$
where $\chi_{i_1,\dots, i_m}:(\dd^*)^m\to \dd^*$
denotes the character  $\lambda_1^{i_1}\cdots \lambda_m^{i_m}$. 
Since the characters of a group (in this case $(\c^*)^m$) are linearly independent we see that 
$$
f_{i_1,\dots, i_m}(x_{m+1},\dots, x_n)\cdot x_1^{i_1}\cdots x_m^{i_m}\in I^{\rm an}+
(x_1,\dots,x_m)^N
$$
holds for every $N$. By Krull's intersection theorem this implies that
$$
f_{i_1,\dots, i_m}(x_{m+1},\dots, x_n)\cdot x_1^{i_1}\cdots x_m^{i_m}\in I^{\rm an}.
$$
We next use translations by $(c_{m+1},\dots, c_m)$ in the  $ x_{m+1},\dots, x_n$ directions to achieve that
$ f_{i_1,\dots, i_m}(x_{m+1}+c_{m+1},\dots, x_n+c_n)$ is nonzero at 
$(x_{m+1},\dots, x_n)=(0,\dots, 0)$. Thus
$$ 
 x_1^{i_1}\cdots x_m^{i_m}\in I^{\rm an} \qtq{provided} f_{i_1,\dots, i_m}(x_{m+1},\dots, x_n)\not\equiv 0.
$$
This shows that $I^{\rm an}$ is generated by monomials in 
$x_1,\dots, x_m$ hence it is toroidal. \qed

\medskip

Note that $(X, \Delta)$ is toric with torus $T$ then  
we need only the $T$-action in the above proof. Thus we have showed the following elementary observation.

\begin{cor} \label{toric.-.toroidal.cor}
Let $(X, \Delta)$ be a smooth toric variety.  
Then an ideal is toric iff it is toroidal.\qed
\end{cor}

Now we come to the key definition, the toroidal hull of an ideal. The existence of the  toroidal hull is a quite elementary observation which is at least implicit in several papers. See, for instance, the notion of the Newton polygon  \cite{kous} and its connections with resolutions \cite{MR2132653} 
or the 
D-balanced  and  well-tuned  ideals discussed in \cite{wlod-res}; see also
 \cite[Sec.3.4]{k-res} for more details on the latter.

\begin{defn-thm} \label{TOR.HULL.EXISTS}
Let $(X, \Delta)$ be an snc pair over a field  and $J\subset \o_X$ an ideal sheaf. There is a unique, smallest
toroidal ideal sheaf  $J^t\supset J$, called the {\it toroidal hull} of $J$.

Furthermore, if $W\subset X$ is a stratum then
$\mult_W J^t=\mult_WJ$. (A stronger version of this property is established in Lemma \ref{torhull.bu.lem}.)
\end{defn-thm}

Proof. As we noted in Paragraph \ref{tor.ideal.say}, specifying $J^t$ is equivalent to
specifying the semigroups $M_Z$ (\ref{tor.ideal.say}.3) and the latter can be done working in an analytic or formal neighborhood of a point
$p_0\in Z^0$ of an open stratum. 

Then the recipe of constructing $J^t$ follows from
the proof of Proposition \ref{cam.=.tor.prop}:
\begin{enumerate}
\item[($*$)] Take all $f=\sum_{i_1,\dots, i_m}f_{i_1,\dots, i_m}(x_{m+1},\dots, x_n)x_1^{i_1}\cdots x_m^{i_m}\in J$ and add the monomial
$x_1^{i_1}\cdots x_m^{i_m}$ to $ J^t$ whenever $f_{i_1,\dots, i_m}\not\equiv 0$.
\end{enumerate}
This also shows that we have not decreased the 
multiplicity along $Z^0$ since
$$
\mult_{p_0}x_1^{i_1}\cdots x_m^{i_m}=
\inf_{p\in Z^0}\mult_p \bigl(f_{i_1,\dots, i_m}\cdot x_1^{i_1}\cdots x_m^{i_m}\bigr)\geq \inf_{p\in Z^0}\mult_p f.\qed
$$

\begin{cor}\label{torres.char} 
Let $(X, \Delta)$ be an snc pair 
and $J\subset \o_X$  an ideal sheaf. Then
 $J$ is toroidally resolved iff
$J^t=\o_X$. \qed
\end{cor}

The following result says that the toroidal hull commutes with
toroidal blow-ups along  toroidally equimultiple centers.

\begin{lem} \label{torhull.bu.lem}
 Assume that $J$ is   toroidally equimultiple along $Z$.
Then
$$
\bigl(\pi^{-1}_*J\bigr)^t=\pi^{-1}_*(J^t).
$$
\end{lem}

Proof. The question is local on $X$ and we can even replace $X$ by its completion $\hat X_x$. Thus we may assume that $(X, \Delta)$ is toric
with torus $T$ acting on $X$. 
Then $J^t=\tsum_{\tau} \tau^*J$ where we sum of all $\tau \in T$. If $J$ is 
toroidally equimultiple along $Z$ with multiplicity $m$ then the same holds for every $\tau^*J$. Thus
$$
\pi^{-1}_*(J^t)=%\pi^{-1}_*\bigl(\tsum_{\tau} \tau^*J\bigr)=
\o_{B_ZX}(mE)\cdot \pi^*\bigl(\tsum_{\tau} \tau^*J\bigr)=
\tsum_{\tau} \bigl(  \o_{B_ZX}(mE)\cdot \tau^*\pi^*J\bigr)=
\bigl(\pi^{-1}_*J\bigr)^t. \qed
$$

The following observations transforms the toroidal resolution problem for $J$ to the usual resolution problem for its  toroidal hull.
Thus the toroidal hull is a variant of the concept of  {\it tuning  an ideal} used in resolution; see \cite[3.54]{k-res}. 

\begin{prop} \label{main.equiv.prop}
Let $(X, \Delta)$ be an snc pair over a field  and $J\subset \o_X$ an ideal sheaf.
There is a natural equivalence between the following sets.
\begin{enumerate}
\item Toroidally  equimultiple blow-up sequences for $J$.
\item Toroidally  equimultiple blow-up sequences for $J^t$.
\item Equimultiple blow-up sequences for $J^t$.
\end{enumerate}
\end{prop} 

Proof. Proposition \ref{TOR.HULL.EXISTS} shows that 
$J$ is toroidally  equimultiple along a stratum $Z$ iff 
$J^t$ is toroidally  equimultiple along $Z$.
A toroidal ideal is toroidally  equimultiple along a stratum $Z$ iff it is
 equimultiple along $Z$. Thus in all 3 settings the blow-ups allowed at the first step are the same. 

Lemma \ref{torhull.bu.lem} guarantees that
this holds for all subsequent steps by induction. \qed

\begin{say}[Resolution of toroidal ideals] \label{main.thm.class.p}
It has been long known that resolution of toric ideal sheaves is a combinatorial questions that is independent of the characteristic
\cite{kkms, amrt, MR1748623, MR1892938}. However, we need a resolution that is obtained by an
equimultiple blow-up sequence. The original toric references that I could find
do not claim this and the methods do not seem to be designed for this purpose.

Resolution of toric and toroidal varieties and ideals using 
equimultiple blow-up sequences is proved in \cite{MR2219845};
see also \cite{MR2959962, MR2999984}.
Note that our setting is quite a bit easier since for us all strata are smooth.
(This is also the reason why we do not need to worry about imperfect fields.)

One should also note that for toroidal ideals an \'etale-local resolution procedure is automatically combinatorial. So, although this is not stated, the resolution method discussed in \cite{wlod-res} and \cite[Chap.3]{k-res} is combinatorial.
Thus it yields the required resolution procedure for  toroidal ideals over any field.
\end{say}

\begin{say}[Proof of Theorem \ref{main.thm.2}]

By Theorem \ref{main.thm.class} (in characteristic $= 0$) and
Paragraph \ref{main.thm.class.p} (in characteristic  $\neq 0$) there is an
equimultiple blow-up sequence 
$$
\bigl(X_n,  \Delta_n, (J^t)_n\bigr)\to \cdots  \to \bigl(X_0, \Delta_0, (J^t)_0\bigr):=
\bigl(X, \Delta, J^t\bigr)
$$
that resolves  $J^t$.
By Proposition \ref{main.equiv.prop} the same sequence
gives a toroidally  equimultiple blow-up sequence for $J$
$$
\bigl(X_n,  \Delta_n, J_n\bigr)\to \cdots  \to \bigl(X_0, \Delta_0, J_0\bigr):=
\bigl(X, \Delta, J\bigr).
$$
By Lemma \ref{torhull.bu.lem} we know that
$(J_n)^t=(J^t)_n$ and the latter is $\o_{X_n}$ by assumption. 
Thus $J_n$ is toroidally resolved by Lemma \ref{torres.char}.
\end{say}

\begin{ack} I thank S.~Keel for posing the question, 
J.~Huh and J.~Tevelev for insightful comments and references. 
Partial  financial support    was provided  by  the NSF under grant number
 DMS-1362960.
\end{ack}

%\bibliography{refs-main/refs}

\def\cprime{$'$} \def\cprime{$'$} \def\cprime{$'$} \def\cprime{$'$}
  \def\cprime{$'$} \def\cprime{$'$} \def\cprime{$'$} \def\dbar{\leavevmode\hbox
  to 0pt{\hskip.2ex \accent"16\hss}d} \def\cprime{$'$} \def\cprime{$'$}
  \def\polhk#1{\setbox0=\hbox{#1}{\ooalign{\hidewidth
  \lower1.5ex\hbox{`}\hidewidth\crcr\unhbox0}}} \def\cprime{$'$}
  \def\cprime{$'$} \def\cprime{$'$} \def\cprime{$'$}
  \def\polhk#1{\setbox0=\hbox{#1}{\ooalign{\hidewidth
  \lower1.5ex\hbox{`}\hidewidth\crcr\unhbox0}}} \def\cdprime{$''$}
  \def\cprime{$'$} \def\cprime{$'$} \def\cprime{$'$} \def\cprime{$'$}
\providecommand{\bysame}{\leavevmode\hbox to3em{\hrulefill}\thinspace}
\providecommand{\MR}{\relax\ifhmode\unskip\space\fi MR }
\providecommand{\MRhref}[2]{%
  \href{http://www.ams.org/mathscinet-getitem?mr=#1}{#2}
}
\providecommand{\href}[2]{#2}

\bigskip

\noindent  Princeton University, Princeton NJ 08544-1000

{\begin{verbatim} kollar@math.princeton.edu\end{verbatim}}

\end{document}